\input amstex
\documentstyle{amsppt}

\NoRunningHeads
\topmatter
\title The Dirichlet-to-Neumann operator for quantum graphs\endtitle
\author Leonid Friedlander \endauthor
\affil University of Arizona\endaffil
\endtopmatter
\document
\head 1. Introduction\endhead
let $G$ be a compact, connected metric graph. It is a one-dimensional simplicial complex. Let
$V$ be the set of its vertices, and let $E$ be the set of its edges; both sets are finite. We assume
that there are no loops. The multiplicity of a vertex is the number of edges that are incident to it;
the fact that an edge $e$ is incidents to a vertex $v$ will be denoted by $v\sim e$ or $e\sim v$.
We allow multiple edges. If one does not like them, one may introduce artificial vertices of
multiplicity $2$ in the middle of multiple edges; that will make no difference in what follows.
Each edge $e$ of $G$ is a segment, and one assigns a length $l(e)>0$ to it. On an edge $e$,
one can introduce an arc length coordinate $s$ that is the distance from a point on the edge
to the initial vertex. It depends on the choice of orientation, so let us fix orientations of edges;
all the constructions below will be independant of this choice. A function $\psi$ on $G$ is a set
of functions $\psi_e(s)$ on edges. Let all functions $\psi_e$ be differentiable on closed edges.
For a vertex $v$ incident to an edge $e$, we define $\psi_e'(v)=-\psi_e'(0)$ if $v$ is the initial
vertex of $e$ and $\psi_e'(v)=\psi_e'(l(e))$ if it is its terminal vertex. Then we define
$$\psi'(v)=\sum_{e\sim v}\psi_e'(v).$$
If $\psi'(v)=0$ then we say that the function $\psi$ satisfies the Neumann (or Kirchhoff) condition at the vertex $v$.
If for all edges $e$ that are incident to a vertex $v$, the values of $\psi_e$ at $v$ are the same
then the function $\psi$ is continuous at the vertex $v$. One can find a good introduction to quantum graphs in [BK].

Let $V_0$ be a subset of $V$; $V_0=\{v_1,\ldots,v_k\}$. It will be called the boundary of $G$. We introduce two operators,
$\Delta_N$ and $\Delta_D$ that act on functions on $G$. Both take a function $\psi$
to $-\psi''$. Here $(\psi'')_e=(\psi_e)''$. Notice that the value of the second derivative of a function is independent of the choice of orientation. The domain of $\Delta_N$ is the
set of all functions $\psi\in H^2(G)$ such that $\psi'(v)=0$ for all $v\in V$. Here
$$H^2(G)=\oplus_{e\in E}H^2(e)\cap C(G)$$
and $H^2(e)$ is the corresponding Sobolev space. The domain of $\Delta_D$ is the space of
$\psi\in H^2(G)$ such that $\psi'(v)=0$ if $v\in V\setminus V_0$ and $\psi(v)=0$ if $v\in V_0$.
It is well known that both operators are self-adjoint, non-negative, and they have discrete spectrum. Let $0<\lambda_1\leq\lambda_2\leq\cdots$ be eigenvalues of $\Delta_D$ and $0=\mu_<\mu_2\leq\mu_3\leq\cdots$ be eigenvalues of $\Delta_N$. As usual, every eigenvalue is counted as many times as its multiplicity is. Eigenvalues of the operators
$\Delta_N$ and $\Delta_D$ are critial values of the quadratic form 
$$\int_G|\psi'(s)|^2ds=\sum_{e\in E}\int_e|\psi_e'(s)|^2ds$$
subject to the constraint $\int_G|\psi|^2ds=1$, with the domains
$$\Cal Q_N=H^1(G)=\oplus_{e\in E}H^1(e)\cap C(G)$$
and
$$\Cal Q_D=\{\psi\in H^1(G):\ \psi(v)=0\quad\text{for}\ v\in V_0\}.$$
The inclusion $\Cal Q_D\subset\Cal Q_N$ and
$$\text{codim}_{\Cal Q_N}\Cal Q_D=k=|V_0|$$
implies trivial inequalities
$$\mu_j\leq\lambda_j\leq\mu_{j+k} \tag 1$$
for every $j$. Let $N_D(\lambda)=\#\{\lambda_j: \lambda_j<\lambda\}$ and
$N_M(\lambda)=\#\{\mu_j: \mu_j<\lambda\}$. Then inequalities (1) are equivalent to
$$0\leq N_N(\lambda)-N_D(\lambda)\leq k \tag 1$'$ $$
for every $\lambda$.

Let $\lambda\not=\lambda_j$ be a real number. Then, for every $x_1,\ldots,x_k$, the problem
$$\psi''+\lambda\psi=0,\quad \psi'(v)=0\quad\text{for}\ v\in V\setminus V_0;\ 
\psi(v_j)=x_j,\ v_j\in V_0$$
is uniquely solvable. Let the function $\psi$ be its solution. Define $\xi_j=\psi'(v_j)$. One gets an operator that maps a point $x=(x_1,\ldots,x_k)\in\Bbb R^n$ to $\xi=(\xi_1,\ldots,\xi_k)\in\Bbb R^n$. We call it the Dirichlet-to-Neumann operator, and we will denote it by $R(\lambda)$. It can be represented by an $k\times k$ matrix. The operator $R(\lambda)$ is self-adjoint because for two functions $\phi$ and $\psi$ that satisfy the Neumann condition at all vertices $v\in V\setminus V_0$ one has
$$\int_G(\phi''\bar\psi-\phi\bar\psi'')ds=\sum_{v_j\in V_0}(\phi'(v_j)\bar\psi(v_j)-\phi(v_j)\bar\psi'(v_j)).$$
If, in addtion $\phi''+\lambda\phi=0$ and $\psi''+\lambda\psi=0$ then
$$\sum_{j=1}^k\phi'(v_j)\bar\psi(v_j)=\sum_{j=1}^k\phi(v_j)\bar\psi'(v_j).\tag 2$$
Let $n_-(\lambda)$ be the number of negative eigenvalues of $R(\lambda)$. Aword-by-word repetition of the argument from [F] (see also [M]) leads to
$$N_N(\lambda)-N_D(\lambda)=n_-(\lambda).\tag 3$$
In a similar way, one can compare spectra of Laplacians with two different Robin conditions. Namely, let $a\geq 0$, and let $\Delta_a$ be the operator $-d^2/ds^2$ on $G$, with the domain that consists of functions $\psi\in H^2(G)$ satisfying the Neumann condition at all vertices outside of $V_0$ and $\psi'(v_j)+a\psi(v_j)=0$ when $v_j\in V_0$. This operator corresponds
to the quadratic form
$$q_a(\psi)=\int_G|\psi'(s)|^2ds+a\sum_{v_j\in V_0}|\psi(v_j)|^2$$
defined on $H^1(G)$. In particular, $\Delta_N=\Delta_0$, and one can interpret $\Delta_D$ as $\Delta_\infty$. Let $\nu_1(a)\leq\nu_2(a)\leq\cdots$ be eigenvalues of $\Delta_a$. Clearly, $\nu_j(a)$ increase with $a$. Let 
$$N_a(\lambda)=\#\{\nu_j(a):\ \nu_j(a)<\lambda\}.$$
Let $a<b$. By $n_{a,b}(\lambda)$ we denote the number of eigenvalues of $R(\lambda)$ that lie in the interval $[-b, -a)$. Then
$$N_a(\lambda)-N_b(\lambda)=n_{a,b}(\lambda).\tag 4$$
The equality (4) is proved in the same way as (3). Let us recall that, outside of the Dirichle spectrum, $R(\lambda)$ is a decreasing operator-valued function of $\lambda$. It means that
$$\frac{d}{d\lambda}(R(\lambda)x,x)<0$$
for every $x\not=0$. In particular, eigenvalues $\sigma_j(\lambda)$ of $R(\lambda)$ are decreasing functions of $\lambda$. They are all positive when $\lambda<0$.  Notice that $\lambda$ is an eigenvalue of $\Delta_a$ if and only if $-a$ is an eigenvalue of $R(\lambda)$.
As $\lambda$ increases, at the moment when some $\sigma_j(\lambda)$ crosses the level $-a$, $N_a(\lambda)$ increases by $1$; in the same way, $N_b(\lambda)$ increases by $1$ at the moment when it crosses the level $-b$. This accounting implies (4).

The question that I treat in this paper is what are possible restriction for $R(\lambda)$. For example, in the case of the Laplacian in a planar domain, $n_-(\lambda)\geq 1$ for $\lambda>0$ ([F]).  The main result of the paper is
\proclaim{Theorem} Let $\lambda>0$, $k\geq 2$, and let $A$ be a symmetric $k\times k$ matrix with real entries. Then there exists a compact, connected metric graph $G$ and $V_0\subset V$, $|V_0|=k$, such that $R(\lambda)=A$.
\endproclaim
This theorem implies that, in general, one can not say anything about mutual order of eigenvalues of the Dirichlet and Neumann Laplacians beyond the trivial estimate (1). The same is true about mutual order of eigenvalues of the Laplacians with two different Robin conditions.\newline
{\it Acknowledgement.} I am grateful to Minh Kha for discussions that stimulated me to ask the question of which symmetric matrices can be realized as matrices of a Dirichlet-to-Neumann operator.

\head 2. Proof of the Theorem\endhead
It is sufficient to prove the theorem for $\lambda=1$. Indeed, let $\tilde G$ be the metric graph that is obtained from $G$ by multiplying lengths of all its edges by $\sqrt{\lambda}$. Then
$R(\lambda)=\sqrt{\lambda}\tilde R(1)$. From this point, we assume $\lambda=1$.

Let $\Cal R_k$, $k\geq 2$, be the set of all symmetric $k\times k$ matrices with real entries
that correspond to the Dirichlet-to-Neumann operator $R=R(1)$ for some connected, compact metric graph with the boundary $V_0$, $|V_0|=k$. Our goal is to prove that $\Cal R_k=\text{Symm}_k$, the space of sll symmetric matrices with real entries. We start with a simple observation: if $A_1$ and $A_2$ belong to $\Cal R_k$ then $A_1+A_1\in \Cal R_k$. Indeed,
let the graph $G_j$ with the boundary $V_0^{(j)}$ has $A_j$ as its Dirichlet-to-Neumann operator, $j=1,2$. Let $V_0^{(j)}=\{v_{1,j},\ldots, v_{k,j}\}$. We form a graph $G$ by gluing $G_1$ and $G_2$ along their boundaries: vertices $v_{i,1}$ and $v_{i,2}$, $i=1,\ldots, k$ are identified. Clearly, the Dirichlet-to-Neumann operator for $G$ is $A_1+A_2$.

First, I will prove the theorem for $k=2$. The Dirichlet-to-Neumann operator for a segment
of length $l$ is
$$R_l=\pmatrix \cot l &-\frac{1}{\sin l}\\
                          -\frac{1}{\sin l}& \cot l
          \endpmatrix .$$
Therefore,
$$\biggl\{\pmatrix a & b\\ b & a\endpmatrix :\ b^2-a^2=1\biggr\}\subset\Cal R_2.$$
 I claim that every point $(a,b)\in \Bbb R^2$ can be represented as the sum of four points
$(a_j, b_j)$, $j=1,2,3,4$, such that $b_j^2-a_j^2=1$. To prove that it is convenient to make a linear change of variables, $x=a+b$, $y=b-a$. Then $x_jy_1=1$. 

First, every point $(x,y)$ with $xy<0$ can be represented as $(x_,y_1)+(x_2,y_2)$, $x_1y_1=x_2y_2=1$. Let, say, $x>0$, $y<0$. Consider the function
$$f(u)=(x+u)\biggl(y+\frac{1}{u}\biggr)=\frac{x}{u}+yu+xy+1.$$
Notice that $f(u)\to +\infty$ when $u\to 0^+$ and $f(u)\to -\infty$ when $u\to +\infty$. Hence $f(u)=1$ for some $u=-x_1>0$, and $(x,y)=(x_1,1/x_1)+(x_2,1/x_2)$ where $x_2=x-x_1$.

Further, every pont $(x,y)\in \Bbb R^2$ can be represented as $(x_1,y_1)+(x_2,y_2)$ with $x_jy_j<0$. If $x\geq 0$, one can take $x_1=x+1$, $y_1<\min\{y,0\}$. Otherwise, one can take $x_1=x-1$, $y_1>\max\{y,0\}$.
This proves that
$$\pmatrix a& b\\ b & a\endpmatrix\in\Cal R_2\tag 5$$
for every $a,b\in \Bbb R$.

To produce matrices with non-equal diagonal elements, we need an operation of concatination. Let $(G, v_1, v_2)$ and $(H, w_1, w_2)$ be connected metric graphs with boundaries $(v_1,v_2)$ and $(w_1,w_2)$, respectively. By $(G\vee H, v_1, w_2)$ I denote the graph that is obtained from the disjoint union of $G$ and $H$ by gluing $v_2$ and $w_1$; the boundary of the new graph is $(v_1, w_2)$.  Let $a\not=0$, and suppose that the Dirichlet-to-Neumann matrix for $G$ is $\text{diag}(a,a)$
be a scalar matrix; let the Dirichlet-to-Newmann operator for $H$ be $0$. Then the Dirichlet-to-Newmann operator of $G\vee H$ is $\text{diag}(a,0)$. Indeed, let $\psi_G$ be the solution
to $\psi_G''+\psi_G=0$ on $G$ such that $\psi_G(v_1)=1$, $\psi_G(v_2)=0$, and let $\psi_H$ be the solution to $\psi_H''+\psi_H=0$ on $H$ such that $\psi_H(w_1)=0$, $\psi_H(w_2)=c$.
Here $c$ is an arbitrary number.
Then $\psi_G'(v_2)=\psi_H'(w_1)$, and
$$\psi(x)=\cases \psi_G(x), & \text{if}\  x\in G\\
                             \psi_H(x), & \text{if}\  x\in H\endcases$$
is a solution to $\psi''+\psi=0$ on $G\vee H$. Clearly, $\psi'(v_1)=a$ and $\psi'(w_2)=0$.
We conclude that
$$\pmatrix a&0\\0&0\endpmatrix\in\Cal R_2$$
for every $a$. Together with (5), it implies $\Cal R_2=\text{Symm}_2$.

Now we pass to the case of an arbitrary $k>2$.  Let $A=(a_{ij})$ be a $k\times k$ symmetric matrix. Take $k$ vertices $v_1,\cdots,v_k$.  For every pair $(i,j)$, $1\leq i<j\leq k$, we define a $2\times 2$ matrix
$$A_{ij}=\pmatrix c_{ij}&a_{ij}\\ a_{ij}&d_{ij}\endpmatrix$$
with
$$c_{ij}=\cases a_{ii}, &\text{if}\  j=i+1\\
                          0,&\text{otherwise} \endcases$$
and
$$d_{ij}=\cases a_{nn},&\text{if}\  i=n-1,\ j=n\\
                           0,&\text{otherwise}\endcases .$$
Let $G_{ij}$ be a connected graph with a two-point boundary, the Dirichlet-to-Neumann operator of which equals $A_{ij}$.We identify the first boundary point of $G_{ij}$ with $v_i$, and the second boundary point with $v_j$. Then $A$ is the Dirichlet-to-Neumann matrix for the graph $G$ that is the union of all $G_{ij}$.
\qed

\enddocument